\documentclass[12pt]{article}
\usepackage{amsmath,amssymb,amsthm}
\usepackage[mathscr]{eucal}
\usepackage{rotate,graphics,epsfig}
\usepackage{color}
\usepackage{hyperref}

\newtheorem{The}{Theorem}
\newtheorem{Exa}[The]{Example}

\theoremstyle{definition}
\newtheorem{Def}[The]{Definition}
\newtheorem{Rem}[The]{Remark}
\numberwithin{equation}{section}
\numberwithin{The}{section}

\newcommand{\be}{\begin{eqnarray}}
\newcommand{\ee}{\end{eqnarray}}

\newcommand{\by}{\begin{eqnarray*}}
\newcommand{\ey}{\end{eqnarray*}}

\newcommand{\bn}{\begin{enumerate}}
\newcommand{\en}{\end{enumerate}}

\newcommand{\bi}{\begin{itemize}}
\newcommand{\ei}{\end{itemize}}

\def\Fbar{{\overline F}}

\def\frac#1#2{{#1 \over #2}}

\def \P {\mathbb{P}}

\def\xx{\boldsymbol{x}}

\def\ww{\boldsymbol{w}}

\begin{document}
\title{Operator Tail Dependence of Copulas}
\author{
Haijun Li\footnote{{\small\texttt{lih@math.wsu.edu}}, Department of Mathematics and Statistics, Washington State University, 
Pullman, WA 99164, U.S.A. }
}
\date{Revised September 2017}
\maketitle

\begin{abstract}
A notion of tail dependence  based on operator regular variation is introduced for copulas, and the standard tail dependence used in the copula literature is included  as a special case. The non-standard tail dependence with marginal power scaling functions having possibly distinct tail indexes is investigated in detail.  
We show that the copulas with operator tail dependence, incorporated  with regularly varying univariate margins, give rise to a rich class of the non-standard multivariate regularly varying distributions. We also show that under some mild conditions, the copula of a non-standard multivariate regularly varying distribution has the standard tail dependence of order 1. 
Some illustrative examples are given.


\end{abstract}

\section{Introduction}
\label{S1}

Tail dependence describes the amount of dependence in the upper tail or lower tail of a multivariate distribution and has been widely used in extreme value analysis and in quantitative risk management \cite{Joe97, mcneil05}. Tail dependence is often studied by using the copula method, which is used to explore scale-invariant features  for a joint distribution. A common theme is to study decay rates of joint tail probabilities of a random vector using same marginal tail scaling functions \cite{Joe97, JLN10}, leading to the so called {\em standard tail dependence}. In this paper, we focus on using different marginal scaling functions in the analysis of copula tail dependence. We introduce the tail dependence for copulas based on {\em operator regular variation} \cite{MS01}, and explore finer tail dependence structures of copulas hidden in the interior that only emerge from possibly distinct marginal tail scalings. 

The copulas with standard tail dependence, incorporated with tail equivalent regularly varying marginal distributions, coincide with standard multivariate regular variation \cite{LS2009, Li2013}. If regularly varying marginal distributions with possibly distinct tail indexes are incorporated, the copulas with standard tail dependence should correspond to non-standard multivariate regularly varying distributions. We make this intuition precise in this paper by showing that the distributions having copulas with operator tail dependence and regularly varying univariate margins constitute a rich class of non-standard multivariate regularly varying distributions; and on the other hand, under some mild conditions, the copula of a non-standard multivariate regular variation has standard tail dependence. By separating margins and scale-invariant dependence structure, copulas provide a flexible tool in multivariate extreme value analysis.

It is worth mentioning that using different marginal tail scalings can reveal hidden tail dependence patterns that are not evident from using standard, tail equivalent marginal scalings. Using the Marshall-Olkin copula, we show that different marginal tail scalings can lead to more natural tail dependence structure. The rigorousness of our operator tail dependence method is based on the general theory of operator regular variation developed in \cite{MS01}. It is evident, such as in finance, that vector data with heavy tails need not have the
same tail index in every direction and that it may be necessary to consider rotated coordinate
systems using operator norming to detect variations in tail behavior \cite{MS99, MS2013}. It is also worth mentioning that intermediate tail dependence in extremes hidden in the interior was observed 
in \cite{LT96, LT98}, and the concept of hidden regular variation in standard form was initiated in \cite{Resnick02}. 
Non-standard multivariate regular variation was addressed in \cite{MS2013} using operator norming, and the application of operator regular variation to geometric extreme value theory and high-dimensional risk analysis can be found in \cite{BE07}.

The rest of this paper is organized as follows. The operator tail dependence of copulas with respect to a matrix of tail indexes is introduced in Section 2, and the main results are proved  in Section 3. Some comments in Section 4 conclude the paper. In what follows, 
two functions $f,g: {\mathbb{R}}\to \mathbb{R}$ are said to be tail equivalent, denoted by  $f(x)\sim g(x)$ as $x\to a$, $a\in \overline{\mathbb{R}}=\mathbb{R}\cup \{+\infty\}$, if $\lim_{x\to a}[f(x)/g(x)]= 1$. A univariate Borel-measurable function $V: \mathbb{R}_+\to \mathbb{R}_+$ is said to be regularly varying at $\infty$ with tail index $\rho\in \mathbb{R}$, denoted by $V\in \mbox{RV}_\rho$, if $V(tx)/V(t)\to x^\rho$ as $t\to \infty$ for any $x>0$. In particular, a function $V\in \mbox{RV}_0$ is said to be slowly varying at $\infty$. See, for example,  \cite{Resnick87, Resnick07} for details on these notions and on the theory of regular variation. For any two vectors ${\boldsymbol {a, b}}\in \mathbb{R}^d$, the sum $\boldsymbol{a+b}$  and vector inequalities such as $\boldsymbol{a\le b}$ are all operated component-wise, and the intervals such as $[{\boldsymbol {a, b}}]$ are the Cartesian products of component-wise intervals.  Without loss of generality, we consider the cone $\mathbb{E}^{(1)}:= [0,\infty]^d\backslash \{\boldsymbol 0\}$ and its interior  $\mathbb{E}^{(2)}:=[0,\infty]^d\backslash \cup_{i=1}^d\{s{\boldsymbol e_i, s\ge 0}\}$, where ${\boldsymbol e}_i$, $i\le i\le d$, denote the standard basis of $\mathbb{R}^d$.

\section{Operator Tail Dependence of Copulas}
\label{operator}

Let $F$ denote a $d$-dimensional distribution function with continuous marginal distributions $F_1, \dots, F_d$. Let $\Fbar_1,\ldots,\Fbar_d$
denote the corresponding marginal survival functions.  The scale-invariant tail behavior of a multivariate distribution $F$ can be studied using its copula, which is invariant under marginal increasing transforms and preserves the scale-invariant dependence structure of the distribution. 

Formally, a copula $C$ is a multivariate distribution with standard uniformly distributed margins on $[0,1]$. Sklar's theorem (see, e.g., Section 1.6 in \cite{Joe97}) states that every multivariate distribution $F$ with margins $F_1,\dots, F_d$ can be written as 
$$F(x_1,\dots, x_d)=C(F_1(x_1),\dots,F_d(x_d)), \ (x_1, \dots, x_d)\in \mathbb{R}^d,$$
for some $d$-dimensional copula $C$. In fact, in the case of continuous margins, $C$ is unique and 
$$C(u_1,\dots,u_d)={F}(F_1^{-1}(u_1), \dots, F_d^{-1}(u_d)), \ (u_1, \dots, u_d)\in [0,1]^d, $$
where $F_i^{-1}(u_i)$ is the quantile function of the $i$-th margin, $1\le i\le d$. 
Let $(U_1,\dots,U_d)$ denote a random vector with distribution $C$ and $U_i, 1\le i\le d$, being uniformly distributed on $[0,1]$. The survival copula $\widehat{C}$ is defined as follows:
\begin{equation}
	\label{dual}
	\widehat{C}(u_1,\dots,u_n)=\P(1-U_1\le u_1,\dots,1-U_n\le u_n)=\overline{C}(1-u_1,\dots,1-u_n)
\end{equation}
where $\overline{C}$ is the joint survival function of $C$. The survival copula $\widehat{C}$ can be used to transform lower tail properties of $(U_1,\dots,U_d)$ into the corresponding upper tail properties of $(1-U_1,\dots,1-U_d)$. 

The upper tail dependence function of $C$ with tail order $\kappa_U$, introduced in \cite{HJ10},
is defined as follows:
\begin{equation}
	b_U({\boldsymbol w}; \kappa_U):= \lim_{u\to 0^+}
	\frac{\overline{C}(1-uw_i, 1\le i\le d)}{u^{\kappa_U}\ell(u)}>0, ~ {\boldsymbol w}=(w_1,\dots,w_d)\in [0,\infty)^d\backslash \{\boldsymbol 0\},\label{U tail}
\end{equation}
provided that the non-zero limit exists for some $\kappa_U\ge 1$ and some function $\ell(\cdot)$ that is slowly varying at $0$; that is, $\ell(s^{-1})\in \mbox{RV}_0$. Similarly, the lower {tail dependence function of $C$ with tail order $\kappa_L$} 
is defined as follows:
\begin{equation}
	b_L({\boldsymbol w}; \kappa_L):= \lim_{u\to 0^+}
	\frac{C(uw_i, 1\le i\le d)}{u^{\kappa_L}\ell(u)}>0, ~ {\boldsymbol w}=(w_1,\dots,w_d)\in [0,\infty)^d\backslash \{\boldsymbol 0\},\label{L tail}
\end{equation}
provided that the limit exists for some $\kappa_L\ge 1$ and some  function $\ell(\cdot)$ that is slowing varying at $0$. Clearly, the tail dependence functions are homogeneous; that is,
\begin{equation}
	\label{scaling}
	b_U(t{\boldsymbol w}; \kappa_U)=t^{\kappa_U}\,b_U({\boldsymbol w}; \kappa_U), ~~b_L(t{\boldsymbol w}; \kappa_L)=t^{\kappa_L}\,b_L({\boldsymbol w}; \kappa_L), 
\end{equation}
for any $t>0$ and ${\boldsymbol w}\in [0,\infty)^d\backslash\{\boldsymbol 0\}$.  It follows from \eqref{scaling} that  tail dependence functions  $b_U(\cdot; \kappa_U)$ and $b_L(\cdot; \kappa_L)$ that are positive at some interior points are positive everywhere in the interior of   $[0,\infty)^d\backslash\{\boldsymbol 0\}$ (also see \cite{JLN10, HJL12}). 
Since the upper tail dependence function of a copula $C$ is the lower tail dependence function of the survival copula $\widehat{C}$ (see \eqref{dual}), we focus on the lower tail dependence. 

The tail dependence function \eqref{U tail} with tail order 1 was introduced and studied in \cite{Jaw06, KKP06} and further studied in \cite{LS2009, NJL08, JLN10}, and in particular, the relation between the tail dependence functions and multivariate regular variation was established in \cite{LS2009}. Various tail dependence parameters used in the copula literature (see, e.g., \cite{Li2009}) can actually be written in terms of the tail dependence functions. The higher order tail dependence functions \eqref{U tail} and \eqref{L tail} were studied in \cite{HJ10}, and the relations between the higher order tail dependence functions and hidden regular variation were established in \cite{HJL12}. The relations between the copula tail dependence and multivariate regular variation were investigated in details in \cite{LW2013, LH14} via copula tail densities. The tail dependencies of copulas have been widely applied to tail risk assessment (see, e.g., \cite{JL11} and the references therein). 

The lower tail dependence function \eqref{L tail} can be written as
\[	b_L({\boldsymbol w}; \kappa_L):= \lim_{u\to 0^+}
\frac{C(u^{1/\kappa_L}w_i, 1\le i\le d)}{u\,\ell(u)}, ~ {\boldsymbol w}=(w_1,\dots,w_d)\in [0,\infty)^d\backslash \{\boldsymbol 0\},
\]
for some function $\ell(\cdot)$ that is slowly varying at $0$. That is, all the margins converge with the same scaling function $u^{1/\kappa_L}$ and common tail index $1/\kappa_L$. In many cases (see \cite{MS99, Resnick07}), margins may converge with scaling functions that have different tail indexes, and in these situations, more subtle normalizations would be needed to reveal finer extremal dependence structure. For example, one can study the asymptotic behavior of $C(u^{\lambda_1}w_1, \dots, u^{\lambda_d}w_d)$ with distinct $\lambda_1, \dots, \lambda_d$, as $u\to 0$, which would lead to {\em non-standard regular variation} \cite{Resnick07, MS2013}. We here adopt a more general approach based on operator regular variation developed in \cite{MS01} and introduce the operator tail dependence of a copula with respect to a diagonalizable matrix.     

Given a $d\times d$ matrix $A$, we define the exponential matrix
\[\exp(A)=\sum_{k=0}^\infty \frac{A^k}{k!}, \ \mbox{where}\ A^0=I\ \mbox{(the $d\times d$ identity matrix)}, 
\]
and the power matrix
\begin{equation}
\label{power matrix}
u^A=\exp(A\log u)= \sum_{k=0}^\infty \frac{A^k(\log u)^k}{k!}, \ \mbox{for}\ u>0.
\end{equation}
Power matrices behave like power functions; for example, for any positive-definite matrix $A$ and any norm $||\cdot||$ on $\mathbb{R}^d$, $||u^A\ww||\to 0$, as $u\to 0$, uniformly on compact subsets of $\ww\in \mathbb{R}^d$. A good summary on properties of exponential and power matrices can be found in \cite{MS01}. 
Note that $u^{-A}=(u^{-1})^A=(u^A)^{-1}$, and thus the limiting results in \cite{MS01} for $t\to \infty$ can be converted to the similar results for $u\to 0$.  

For any positive-definite matrix $A$, it follows from the eigen-decomposition that 
\[A = O^{-1}
\begin{pmatrix}
\lambda_1& \cdots & 0\\
 \vdots & \ddots & \vdots \\
 0 & \cdots & \lambda_d
\end{pmatrix}
O,
\]
where $O$ is an orthogonal matrix and eigenvalues $\lambda_1, \dots, \lambda_d$ are all positive. Then it is easy to see that
\begin{equation}
u^A=O^{-1}
\begin{pmatrix}
u^{\lambda_1} & \cdots & 0\\
\vdots & \ddots & \vdots \\
0 & \cdots & u^{\lambda_d}
\end{pmatrix}
O,\ \ \mbox{for}\ u>0. 
\label{dec}
\end{equation}
Let ${\boldsymbol W}^A = \{\ww\in [0,\infty)^d\backslash \{\boldsymbol 0\}: u^{A}\ww\ge {\boldsymbol 0} \ \mbox{for all sufficiently small } u\}$.

\begin{Def}
	\label{operator a}
	Let $C$ be a $d$-dimensional copula of ${\boldsymbol U}=(U_1, \dots, U_d)$ with standard uniform margins $U_i$'s, and let $A$ be a $d\times d$ positive-definite matrix. The lower and upper operator exponent functions of $C$ with respect to matrix index $A$ are defined as follows, for a vector ${\boldsymbol w}=(w_1,\dots,w_d)\in [0,\infty)^d\backslash \{\boldsymbol 0\}$, 
	\begin{equation}
	a_L({\boldsymbol w}; A, C) := \left\{
	\begin{array}{ll}
	\lim_{u\to 0^+}
	\frac{\P\big(\boldsymbol U\in (u^{A}\ww, {\boldsymbol 1}]^c\big)}{u\,\ell(u)}, & \text{for } \ww \in {\boldsymbol W}^A\\
	0, & \text{otherwise;} \label{operator a_L}
	\end{array}
	\right.
	\end{equation}
	\begin{equation}
	a_U({\boldsymbol w}; A, C) := a_L({\boldsymbol w}; A, \widehat{C})= \left\{
	\begin{array}{ll}
	\lim_{u\to 0^+}
	\frac{\P\big(\boldsymbol U\in [{\boldsymbol 0}, {\boldsymbol 1}-u^{A}\ww]^c\big)}{u\,\ell(u)},& \text{for } \ww \in {\boldsymbol W}^A\\
	0, & \text{otherwise.}\label{operator a_U}
		\end{array}
		\right.
	\end{equation}
	provided that the limits exist for some function $\ell(\cdot)$ that is slowly varying at $0$. Here, and in the sequel, the vectors ${\boldsymbol 0}$ and ${\boldsymbol 1}$ denote the vectors of zeros and 1's respectively, and $S^c$ denotes the complement of subset $S$ within the support of the underlying measure.  
\end{Def}

\begin{Rem}
	\begin{enumerate}
		\item If $A$ is a diagonal matrix with positive diagonal entries $\lambda_1, \dots, \lambda_d$, then \eqref{operator a_L} and \eqref{operator a_U} reduce, respectively, to
		\be
		a_L(\ww; A,C)&=&\lim_{u\to 0^+}\frac{\P(U_i\le u^{\lambda_i}w_i, \mbox{for some}\ 1\le i\le d)}{u\ell(u)}\label{non-standard a_L}\\
		a_U(\ww; A,C)&=&\lim_{u\to 0^+}\frac{\P(U_i> 1-u^{\lambda_i}w_i, \mbox{for some}\ 1\le i\le d)}{u\ell(u)}.\label{non-standard a_U}
		\ee 
		The case that $\lambda_1=\dots = \lambda_d$ is discussed in details in \cite{HJ10, HJL12}. 
		\item The operator exponent functions \eqref{operator a_L} and \eqref{operator a_U} may not exist at some points in the cone $\mathbb{E}^{(1)}$. For example, if $A$ is a diagonal matrix with positive diagonal entries that are all strictly less than 1, then the exponent functions \eqref{non-standard a_L} and \eqref{non-standard a_U} do not exist on the  axes $\{s{\boldsymbol e_i, s\ge 0}\}$, $1\le i\le d$, and in these cases, the exponent functions defined on the subcone $\mathbb{E}^{(2)}$ only measure hidden dependence in the interior $\mathbb{E}^{(2)}$.
	\end{enumerate}
\end{Rem}

\begin{Rem}
	\begin{enumerate}
		\item Instead of using $u^A$ in \eqref{operator a_L} and \eqref{operator a_U}, one can use a general scaling function that is operator regularly varying with matrix index $A$ (see \cite{MS01}). Using $u^A$ with a positive-definite matrix $A$ suits our purpose of studying higher order non-standard tail dependence of copulas. It is also worth mentioning that the slowly varying function $\ell(\cdot)$ is not unique but such functions are tail equivalent. 
		\item The scaling properties of operator exponent functions can be stated, in terms of the matrix $A$, as follows.
		\begin{equation}
		\label{operator scaling}
		a_L(t^{A}{\boldsymbol w}; A, C)=t\,a_L({\boldsymbol w}; A, C),\ \ a_U(t^{A}{\boldsymbol w}; A, C)=t\,a_U({\boldsymbol w}; A, C).
		\end{equation}
		This is due to, for the lower case, the fact that 
		\[a_L(t^{A}{\boldsymbol w}; A, C)=\lim_{u\to 0^+}
		\frac{\P(\boldsymbol U\in (u^{A}t^A\ww, {\boldsymbol 1}]^c)}{u\,\ell(u)}=t\lim_{u\to 0^+}
		\frac{\P(\boldsymbol U\in ((ut)^A\ww, {\boldsymbol 1}]^c)}{tu\,\ell(tu)}=t\,a_L({\boldsymbol w}; A, C).
		\]
		The upper case is similar. When $A$ is diagonal, \eqref{operator scaling} is often called ``weighted homogeneity'' or ``quasihomogeneity'' in the literature (see \cite{AGV12} and references therein). 
	\end{enumerate}
\end{Rem}

The exponent function \eqref{operator a_L} can be used to estimate tail probabilities  $\P(\boldsymbol U\in u^AB)$ for a  Borel set $B \subset [0,1]^d$ that is bounded away from $\boldsymbol 1$, where $u^AB:= \{u^A\ww: \ww\in B\}$. In particular, orthant sets such as $B=[{\boldsymbol 0}, \ww]$ are frequently encountered in tail risk analysis \cite{BE07, JL11}, and so the tail dependence functions become particularly useful.  

\begin{Def}
	\label{operator tail}
	Let $C$ be a $d$-dimensional copula and $A$ be a $d\times d$ positive-definite matrix. The lower and upper tail dependence functions of $C$ with respect to matrix index $A$ are defined as follows, for a vector ${\boldsymbol w}=(w_1,\dots,w_d)\in [0,\infty)^d\backslash \{\boldsymbol 0\}$, 
	\be
	b_L({\boldsymbol w}; A, C) &:=& \lim_{u\to 0^+}
	\frac{C(u^{A}\ww)}{u\,\ell(u)},\ \text{for } \ww \in {\boldsymbol W}^A\label{operator lower}\\
		b_U({\boldsymbol w}; A, C) &:=& \lim_{u\to 0^+}
		\frac{\overline{C}({\boldsymbol 1}-u^{A}\ww)}{u\,\ell(u)},\ \text{for } \ww \in {\boldsymbol W}^A\label{operator upper}
	\ee
	for some function $\ell(\cdot)$ that is slowly varying at $0$, and zero otherwise. 
\end{Def}

\begin{Rem}
	\begin{enumerate}
		\item If $A$ is a diagonal matrix with identical, positive diagonal entries, then \eqref{operator lower} and \eqref{operator upper} reduce to \eqref{L tail} and \eqref{U tail} respectively. 
		\item Most of the extremal dependence measures used in practice are the special cases of \eqref{operator lower} or \eqref{operator upper}. Unfortunately, the existences of tail dependence functions \eqref{operator lower} and \eqref{operator upper} do not guarantee the existences of tail dependence functions of lower dimensional margins \cite{HJL12}. As was shown in Lemma 6.1 of \cite{Resnick07} (also see \cite{HJL12}),   the existence of exponent functions \eqref{operator a_L} and \eqref{operator a_U} is equivalent to the vague convergence, which ensures the existence of all marginal exponent functions and the existence of tail dependence functions \eqref{operator lower} and \eqref{operator upper}.  
	\end{enumerate}
\end{Rem}

We illustrate these functions \eqref{operator lower} and \eqref{operator upper} using the bivariate Marshall-Olkin copula. 
Let $\{E_1, E_2, E_{12}\}$ be a sequence of
independent,  exponentially distributed random variables, with
$E_S$ having mean $1/\lambda_S$, $\lambda_S> 0$, $\emptyset\ne S\subseteq \{1, 2\}$. Let
\begin{equation}
\label{MO} 
T_1 = \min\{E_1, E_{12}\}, \ T_2 = \min\{E_2, E_{12}\}
\end{equation}
The joint distribution of ${\boldsymbol T} =(T_1, T_2)$ is
called the bivariate Marshall-Olkin exponential distribution with parameters
$\{\lambda_1, \lambda_2, \lambda_{12}\}$ (see \cite{MO67}). In the reliability context, $T_1, T_2$ can be viewed
as the lifetimes of two components operating in a random shock
environment where a fatal shock governed by Poisson process $\{N_S(t), t\ge 0\}$
with
rate $\lambda_S$ destroys all the components with indexes in
$S\subseteq \{1, 2\}$ simultaneously. 
It follows from \eqref{MO} that the survival function of ${\boldsymbol T}$ can be written as
\be
\overline{F}(t_1, t_2) &=& \P(T_1>t_1, T_2>t_2)= \exp\left [-\lambda_1t_1 -\lambda_2t_2 - \lambda_{12}\max\{t_1,t_2\}
\right ].\label{MO survival}
\ee
Its survival copula is derived as follows
\begin{equation}
\label{survival copula}
\widehat{C}(u_1, u_2) =  u_1u_2\min\{u_1^{-\alpha_1^{12}}, u_2^{-\alpha_2^{12}}\},\ \ 0\le u_1, u_2\le 1, 
\end{equation}
	where $\alpha_i^{12}=\lambda_{12}/(\lambda_i+\lambda_{12})$, $\lambda_{12}>0$, $\lambda_i>0$, $i=1, 2$. 
The general Marshall-Olkin survival copula in high dimension is derived in \cite{Li2008}. 
The upper tail dependence function of the Marshall-Olkin survival copula $\widehat{C}$ with tail order $\kappa_U=1$ is also derived in \cite{Li2008}. It is shown in \cite{Li2008} that the lower tail dependence function of  $\widehat{C}$ with tail order $\kappa_L=1$ is zero. To reveal higher order regular variation properties for lower tails of $\widehat{C}$ hidden within the interior of $(0,\infty)^2$, we use the operator regular variation. 
\begin{Exa}\rm
	\label{2}
	The survival copula of a bivariate Marshall-Olkin vector ${\boldsymbol T} = (T_1, T_2)$ is given by \eqref{survival copula}.  Consider
	\by
	\widehat{C}(uw_1, uw_2) &=&  u^2w_1w_2\min\{u^{-\alpha_1^{12}}w_1^{-\alpha_1^{12}}, u^{-\alpha_2^{12}}w_2^{-\alpha_2^{12}}\}\\
	&=& u^{2-\alpha_1^{12}}w_1w_2\min\{w_1^{-\alpha_1^{12}}, u^{\alpha_1^{12}-\alpha_2^{12}}w_2^{-\alpha_2^{12}}\}\\
	&\sim &
	\left\{
	\begin{array}{ll}
		u^{2-\alpha_1^{12}}w_1^{1-\alpha_1^{12}}w_2 & \ \ \ \mbox{if $\alpha_1^{12}< \alpha_2^{12}$}\\
		u^{2-\alpha_1^{12}}w_1w_2\min\{w_1^{-\alpha_1^{12}}, w_2^{-\alpha_2^{12}}\} & \ \ \ \mbox{if $\alpha_1^{12}= \alpha_2^{12}$}\\
		u^{2-\alpha_2^{12}}w_1w_2^{1-\alpha_2^{12}} & \ \ \ \mbox{if $\alpha_1^{12}>\alpha_2^{12}$,}
	\end{array}
	\right.
	\ey
as $u\to 0$. Therefore,
\begin{equation}
\label{bi-1}
b_L(\ww; A_1, \widehat{C})=\lim_{u\to 0^+}\frac{\widehat{C}(uw_1, uw_2)}{u^{2-\min\{\alpha_1^{12},\alpha_2^{12}\}}}
=	\left\{
\begin{array}{ll}
w_1^{1-\alpha_1^{12}}w_2 & \ \ \ \mbox{if $\alpha_1^{12}< \alpha_2^{12}$}\\
w_1w_2\min\{w_1^{-\alpha_1^{12}}, w_2^{-\alpha_2^{12}}\} & \ \ \ \mbox{if $\alpha_1^{12}= \alpha_2^{12}$}\\
w_1w_2^{1-\alpha_2^{12}} & \ \ \ \mbox{if $\alpha_1^{12}>\alpha_2^{12}$.}
\end{array}
\right.
\end{equation}
That is, the lower tail order $\kappa_L=2-\min\{\alpha_1^{12},\alpha_2^{12}\}$ and the lower tail dependence function of $\widehat{C}$ is given by \eqref{bi-1}. To rephrase this in terms of operator regular variation, the Marshall-Olkin survival copula has a lower tail dependence function with matrix index 
\[A_1=\begin{pmatrix}
\lambda &  0\\
0 & \lambda 
\end{pmatrix},
\ \ \mbox{where}\ \lambda = \frac{1}{\kappa_L}=\frac{1}{2-\min\{\alpha_1^{12},\alpha_2^{12}\}}.
\]

We now use a different power matrix scaling, that leads to uncovering a natural extremal dependence structure. Consider
\by
\widehat{C}(u^{1/\alpha_1^{12}}w_1, u^{1/\alpha_2^{12}}w_2)&=&u^{1/\alpha_1^{12}}u^{1/\alpha_2^{12}}w_1w_2\min\{u^{-1}w_1^{-\alpha_1^{12}}, u^{-1}w_2^{-\alpha_2^{12}}\}\\
&=& u^{\frac{\lambda_1+\lambda_2+\lambda_{12}}{\lambda_{12}}}w_1w_2\min\{w_1^{-\alpha_1^{12}}, w_2^{-\alpha_2^{12}}\}.
\ey
Let $\beta_1=\frac{\lambda_1+\lambda_{12}}{\lambda_1+\lambda_2+\lambda_{12}}$ and $\beta_2=\frac{\lambda_2+\lambda_{12}}{\lambda_1+\lambda_2+\lambda_{12}}$, and we have,
\begin{equation}
\label{bi-2}
b_L(\ww; A_2, \widehat{C})=\lim_{u\to 0^+}\frac{\widehat{C}(u^{\beta_1}w_1, u^{\beta_2}w_2)}{u}=w_1w_2\min\{w_1^{-\alpha_1^{12}}, w_2^{-\alpha_2^{12}}\}, \ \ w_1>0, w_2>0.
\end{equation}
That is, the Marshall-Olkin survival copula has a lower tail dependence function with matrix index 
\begin{equation}
\label{bi-3}
A_2=\begin{pmatrix}
\beta_1 &  0\\
0 & \beta_2 
\end{pmatrix}.
\end{equation}
In contrast to \eqref{bi-1}, \eqref{bi-2} resembles the same structure as that of the Marshall-Olkin survival copula. Observe that both $\beta_1$ and $\beta_2$ are marginal parameters and $\alpha_1^{12}$ and $\alpha_2^{12}$ are the parameters that describe the (asymmetric) dependence structure. 
	\hfill $\Box$
\end{Exa}

The Marshall-Olkin distribution has no upper tail dependence of first order \cite{Li2008}. Furthermore, Example \ref{2} shows that the higher order upper tail dependence of the Marshall-Olkin distribution exhibits different structures with different marginal scalings.

\section{Non-Standard Regular Variation and Tail Dependence of Copulas}
\label{NRV}

In this section, we show that the copulas with operator exponent functions, incorporated with regularly varying univariate margins, constitute a rich class of non-standard multivariate regularly varying distributions. On the other hand, the copula of a non-standard regularly varying distribution has a (standard) tail dependence of first order under some mild regularity conditions.

Consider again a non-negative random vector $\boldsymbol X=(X_1, \dots, X_d)$ with joint distribution $F$ having continuous marginal distributions $F_1, \dots, F_d$. The component index set is denoted by $D=\{1, \dots, d\}$. 
\begin{Def}
	\label{non-standrad RV}
	 The random vector $\boldsymbol X$ is said to be non-standard regularly varying if there exists a Radon measure $\mu(\cdot)$
	(i.e., finite on compact sets), called the {\em intensity measure},  on $\mathbb{E}^{(1)}=\mathbb{\overline R}^d_+\backslash\{\boldsymbol{0}\}$, a regularly varying function $R(t)\in \mbox{RV}_{-\beta}$, $\beta>0$, and positive marginal scaling indexes $\gamma_1, \dots, \gamma_d$ such that  
	\begin{equation}
	\lim_{t\to \infty}\frac{\P(X_i> t^{\gamma_i}x_i,\ \exists i\in D)}{R(t)}=\mu([{\boldsymbol 0}, \xx]^c)
	\label{non-stand rv}
	\end{equation}
	for every continuity point $\xx=(x_1, \dots, x_d)\in \mathbb{\overline R}^d_+\backslash\{\boldsymbol{0}\}$ of the function $\mu([{\boldsymbol 0}, \cdot]^c)$. 
\end{Def}

\begin{Rem}
	\begin{enumerate}
		\item If we let
		\[E=\begin{pmatrix}
		\gamma_1& \cdots & 0\\
		\vdots & \ddots & \vdots \\
		0 & \cdots & \gamma_d
		\end{pmatrix}.
		\]
		then \eqref{non-stand rv} can be rephrased as
		\begin{equation}
		\label{operator non-stan}
		\lim_{t\to \infty}\frac{\P\big(\boldsymbol X\in [{\boldsymbol 0}, t^E\xx]^c\big)}{R(t)}=\mu([{\boldsymbol 0}, \xx]^c),
		\end{equation}
		for every continuity point $\xx=(x_1, \dots, x_d)\in \mathbb{\overline R}^d_+\backslash\{\boldsymbol{0}\}$. In general, an operator regularly varying scaling function can be used in place of $t^E$ \cite{MS01}. 
		\item It is known that \eqref{non-stand rv} is equivalent to the vague convergence of Radon measures (see \cite{Resnick07}, Lemma 6.1). 
		The intensity measure $\mu(\cdot)$ also enjoys the following scaling property:
		\[\mu(s^EB)=s^{-\beta}\mu(B),\ s>0, 
		\]
		for any Borel subset $B\subset \mathbb{\overline R}^d_+\backslash\{\boldsymbol{0}\}$, with $\mu(\partial B)=0$, that is bounded away from the origin. 
		\item The case that $\gamma_1=\dots = \gamma_d$ is called the standard regular variation \cite{Resnick07}. In contrast to the standard regular variation, the non-standard regular variation allows possibly distinct tail indexes among the margins. 
	\end{enumerate}
\end{Rem}

\begin{Rem}
	Non-standard multivariate regular variation can be defined similarly on $\mathbb{\overline R}^d\backslash\{\boldsymbol{0}\}$ \cite{MS2013}, with a strong motivation from financial data analysis. Currency exchange rates of Deutsche Mark and Yen versus the US Dollar were observed with the preponderance of probability mass lying near the diagonal
	line with slope +1, and the data were analyzed in \cite{NPM01} by fitting a multivariate stable model to the data, assuming that a uniform tail
	thickness in every radial direction. The same data were analyzed in \cite{MS99} using operator norming and found that the tail thickness varies significantly
	with direction. Specifically, the approach based on eigenvectors of the sample covariance matrix yields new coordinates on $\mathbb{\overline R}^2\backslash\{\boldsymbol{0}\}$ that determine completely the tail behavior and tail thickness estimates. 
\end{Rem}

Let $C$ denote the copula of distribution $F$ defined on $\mathbb{R}^d_+$, and  $\Fbar_1, \dots, \Fbar_d$ denote the univariate,  continuous marginal survival functions.

\begin{The}\rm 
	\label{TDF-MRV}
	Suppose that the copula $C$ has the upper operator exponent function $a_U(\cdot; A, C)$ where 
	\[A=\begin{pmatrix}
	\lambda_1& \cdots & 0\\
	\vdots & \ddots & \vdots \\
	0 & \cdots & \lambda_d
	\end{pmatrix}.
	\]
	If $a_U(\cdot; A, C)$ is continuous, and if, for all  $1\le i\le d$, $\Fbar_i(t)=t^{-\alpha_i}L_i(t)\in \mbox{RV}_{-\alpha_i}$, $\alpha_i>0$, satisfy that $L_i(t)\to l_i>0$,   
	then the distribution $F(x_1, \dots, x_d)=C(F_1(x_1), \dots, F_d(x_d))$ is non-standard regularly varying in the sense of \eqref{non-stand rv}, with $\beta = \alpha_1$ and  marginal scaling indexes $\gamma_i = \lambda_i\alpha_1/\alpha_i$, $i=1, \dots, d$. The intensity measure in \eqref{non-stand rv} is given by
	\[\mu([{\boldsymbol 0}, \ww]^c):=a_U((w_1^{-\alpha_1}r_1, \dots, w_d^{-\alpha_d}r_d); A, {C}),\ r_i=l_i/l_1^{\lambda_i},\ i=1, \dots, d. 
	\] 
\end{The}

\noindent {\bf Proof.} Observe first that $a_L(\cdot; A, \widehat{C})=a_U(\cdot; A, C)$. 
Consider
\by
\P\big(X_i>t^{\lambda_i\alpha_1/\alpha_i}w_i,\ \exists i\in D\big)&=& \P\big(\overline{F}_i(X_i)\le \overline{F}_i(t^{\lambda_i\alpha_1/\alpha_i}w_i), \ \exists i\in D\big)\\
&=& \P\big(\overline{F}_i(X_i)\le
t^{-\lambda_i\alpha_1}w_i^{-\alpha_i}L_i(t^{\lambda_i\alpha_1/\alpha_i}w_i), \ \exists i\in D\big)\\
&=& \P\big(\overline{F}_i(X_i)\le
(t^{-\alpha_1}L_1(t))^{\lambda_i}w_i^{-\alpha_i}R_i(t), \ \exists i\in D\big)\\
&=& \P\big(\overline{F}_i(X_i)\le
\overline{F}_1^{\lambda_i}(t)w_i^{-\alpha_i}R_i(t), \ \exists i\in D\big)
\ey
where $R_i(t)=L_i(t^{\lambda_i\alpha_1/\alpha_i}w_i)/L_1^{\lambda_i}(t)$, $1\le i\le d$. Since $R_i(t)\to r_i:=l_i/l_1^{\lambda_i}>0$, as $t\to \infty$, we have for any small
$\epsilon>0$, when $t$ is sufficiently large,
\[r_i-\epsilon < R_i(t) < r_i+\epsilon, ~\mbox{for all}~ 1\le i\le d.
\]
That is, when $t$ is sufficiently large,
\[\P\big(\overline{F}_i(X_i)\le\overline{F}_1^{\lambda_i}(t)w_i^{-\alpha_i}(r_i-\epsilon), \ \exists i\in D\big)
\]
\[\le \P\big(X_i>t^{\lambda_i\alpha_1/\alpha_i}w_i, \ \exists i\in D\big)\le \P\big(\overline{F}_i(X_i)\le
\overline{F}_1^{\lambda_i}(t)w_i^{-\alpha_i}(r_i+\epsilon), \ \exists i\in D\big).
\]
Let $u=\overline{F}_1(t)$. Since $a_L(\cdot; A, \widehat{C})$ exists, the following two limits are well-defined:
\[\lim_{u\to 0^+}\frac{\P\big(\overline{F}_i(X_i)\le u^{\lambda_i}w_i^{-\alpha_i}(r_i-\epsilon), \ \exists i\in D\big)}{u\ell(u)}=a_L((w_1^{-\alpha_1}(r_1-\epsilon), \dots, w_d^{-\alpha_d}(r_d-\epsilon)); A, \widehat{C})
\]
\[\lim_{u\to 0^+}\frac{\P\big(\overline{F}_i(X_i)\le u^{\lambda_i}w_i^{-\alpha_i}(r_i+\epsilon), \ \exists i\in D\big)}{u\ell(u)}=a_L((w_1^{-\alpha_1}(r_1+\epsilon), \dots, w_d^{-\alpha_d}(r_d+\epsilon)); A, \widehat{C}).
\]
Therefore, we have, for any small
$\epsilon>0$,
\[a_L((w_1^{-\alpha_1}(r_1-\epsilon), \dots, w_d^{-\alpha_d}(r_d-\epsilon)); A, \widehat{C})\le 
\liminf_{t\to \infty}\frac{\P\big(X_i>t^{\lambda_i\alpha_1/\alpha_i}w_i, \ \exists i\in D\big)}{\overline{F}_1(t)\ell(\overline{F}_1(t))}\]
\[
\le \limsup_{t\to \infty}\frac{\P\big(X_i>t^{\lambda_i\alpha_1/\alpha_i}w_i, \ \exists i\in D\big)}{\overline{F}_1(t)\ell(\overline{F}_1(t))}
\le a_L((w_1^{-\alpha_1}(r_1+\epsilon), \dots, w_d^{-\alpha_d}(r_d+\epsilon)); A, \widehat{C}). 
\]
Letting $\epsilon\to 0$, these inequalities and the continuity of $a_U(\cdot; A, C)$ imply that the following limit exists 
\[\lim_{t\to \infty}\frac{\P\big(X_i>t^{\lambda_i\alpha_1/\alpha_i}w_i, \ \exists i\in D\big)}{\overline{F}_1(t)\ell(\overline{F}_1(t))}
= a_L((w_1^{-\alpha_1}r_1, \dots, w_d^{-\alpha_d}r_d); A, \widehat{C})
\]
where $\overline{F}_1(t)\ell(\overline{F}_1(t))\in \mbox{RV}_{-\alpha_1}$. That is, $(X_1, \dots, X_d)$ is non-standard regularly varying with the intensity measure given by $\mu([{\boldsymbol 0}, \ww]^c)= a_L((w_1^{-\alpha_1}r_1, \dots, w_d^{-\alpha_d}r_d); A, \widehat{C})$ and marginal scaling indexes $\gamma_i=\lambda_i\alpha_1/\alpha_i$, $i=1, \dots, d$. 
 \hfill
 $\Box$

\begin{Rem}
	\begin{enumerate}
		\item Multivariate regular variation in the cone $\mathbb{E}^{(1)}$, where the matrix $A$ is the identity  matrix, has been discussed in the literature; see, for example, \cite{Resnick07, LS2009, JLN10} for details. The exponent functions defined on the subcone $\mathbb{E}^{(2)}$ with equal tail indexes $\lambda_i$s describe the extremal dependence hidden in the interior $\mathbb{E}^{(2)}$ \cite{Resnick07, HJL12}.
		\item As illustrated in the proof of this theorem,  the scaling function $R(t)$ and marginal scaling indexes $\gamma_1, \dots, \gamma_d$ in \eqref{non-stand rv} may not be unique for a given intensity measure $\mu(\cdot)$. However, these parameters $\beta$ and $\gamma_i$, $1\le i\le d$, are tail equivalent at infinity through monotone increasing homomorphisms, so that marginal power scalings with tail indexes $\lambda_i/\alpha_i$, $1\le i\le d$, are unique. Statistical tests, such as the one developed in \cite{MS2013}, should be always performed first to see whether or not operator norming is appropriate.  
	\end{enumerate}
\end{Rem}

\begin{Exa}\rm 
	\label{3}
	Consider a bivariate copula $C$ with the following survival (Marshall-Olkin)  copula (see Example \ref{2}):
	\[\widehat{C}(u_1, u_2) =  u_1u_2\min\{u_1^{-\alpha_1^{12}}, u_2^{-\alpha_2^{12}}\},\]
		where $\alpha_i^{12}=\lambda_{12}/(\lambda_i+\lambda_{12})$, $\lambda_{12}>0$, $\lambda_i>0$, $i=1, 2$. 
	That is, 
	\[C(u_1, u_2)=(u_1+u_2-1)+(1-u_1)(1-u_2)\min\{(1-u_1)^{-\alpha_1^{12}},(1-u_2)^{-\alpha_2^{12}}\}, \ \ 0\le u_1, u_2\le 1. 
	\]
	It follows from \eqref{bi-2} that
	\[b_U(\ww; A, C)=b_L(\ww; A, \widehat{C})= w_1w_2\min\{w_1^{-\alpha_1^{12}}, w_2^{-\alpha_2^{12}}\}, 
	\] 
	where the $2\times 2$ matrix $A$ is given by \eqref{bi-3}, 
	 and so its upper operator tail dependence  function with respect to $A$ is continuous on $\mathbb{E}^{(2)}$. Let $F$ be a bivariate distribution with copula $C$ and Pareto marginal distributions, which are given by
	\[F_i(t)=1-(1+t)^{-\alpha_i},\ t\ge 0,\ i=1, 2.
	\]
	Observe that $\overline{F}_i(t)=t^{-\alpha_i}L_i(t)\in \mbox{RV}_{-\alpha_i}$ where $L_i(t)=t^{\alpha_i}/(1+t)^{\alpha_i}\to 1$ as $t\to \infty$. 
	By Theorem \ref{TDF-MRV}, the distribution $F(x_1, x_2)=C(F_1(x_1), F_2(x_2))$ has non-standard regular variation with intensity measure on $\mathbb{E}^{(2)}$
	\[\mu\big((w_1, \infty]\times (w_2, \infty])\big)= w_1^{-\alpha_1}w_2^{-\alpha_2}\min\{w_1^{\alpha_1^{12}\alpha_1}, w_2^{\alpha_2^{12}\alpha_2}\}, \ \ w_1> 0, w_2> 0,
	\]
	where $\alpha_i^{12}=\lambda_{12}/(\lambda_i+\lambda_{12})$ for $i=1, 2$. Note that the distribution is upper tail independent with univariate marginal survival functions $\overline{F}_1(t)$ and $\overline{F}_2(t)$, as $t\to \infty$, and that the intensity measure $\mu(\cdot)$ describes the dependence hidden within the subcone $\mathbb{E}^{(2)}$. 
	 \hfill
	 $\Box$
\end{Exa}

 In Theorem \ref{TDF-MRV}, the assumption that $L_i(t)\to l_i>0$ for $1\le i\le d$ is a technical condition that makes the proof easier. This assumption is mild and many multivariate distributions used in practice satisfy this assumption (see \cite{Joe97, mcneil05}). For general slowly varying functions $L_i(t)$, $1\le i\le d$, one has to modify the definition \eqref{non-stand rv} using a general operator regularly varying scaling function in place of $t^E$; see \cite{MS01} for details on a general theory of the operator regular variation.

 We continue to assume in the rest of this paper that the marginal survival functions $\Fbar_i(t)=t^{-\alpha_i}L_i(t)$, $1\le i\le d$, where the slowly varying function $L_i(t)\to l_i>0$, as $t\to \infty$.
 This implies that for any $i$ and $j$, 
 \begin{equation}
 \label{marginal assumption}
 \frac{\Fbar_i(t^{1/\alpha_i})}{\Fbar_j(t^{1/\alpha_j})}\to \frac{l_i}{l_j}>0.
 \end{equation}
 Suppose that the scaling function $R(t)$ in \eqref{non-stand rv} is written as $R(t)=t^{-\beta}L(t)\in \mbox{RV}_{-\beta}$. Let 
  \begin{equation}
  \label{matrix index A}
  A=\begin{pmatrix}
  \alpha_1\gamma_1/\beta & \cdots & 0\\
  \vdots & \ddots & \vdots \\
  0 & \cdots & \alpha_d\gamma_d/\beta
  \end{pmatrix},
  \end{equation}
  where $\gamma_1, \dots, \gamma_d$ are marginal scaling indexes in \eqref{non-stand rv}.

\begin{The}\rm
	\label{MRV-TDF} Consider a non-negative random vector ${\boldsymbol X}=(X_1,
	\dots, X_d)$ with non-standard regularly varying distribution  $F$ and continuous margins $F_1, \dots,
	F_d$, satisfying that \eqref{non-stand rv} and \eqref{marginal assumption}. 
	Let $C$ and $\mu(\cdot)$ denote, respectively, the copula and the
	intensity measure of $F$. If the intensity measure is orthant-continuous, that is, $\mu([{\boldsymbol 0}, \ww]^c)$ is continuous in $\ww$, then  
	 the upper operator exponent function $a_U(\,\cdot\,; A, C)$ exists
	and
	 \[a_U(\ww; A, C)=\mu\Big(\big(\prod_{i=1}^d[0,w_i^{-1/\alpha_i}l_i^{1/\alpha_i}]\big)^c\Big)
	 \] 
	 for $\ww=(w_1, \dots, w_d)\in \mathbb{R}^d_+\backslash\{\boldsymbol{0}\}$, where the matrix $A$ is given by \eqref{matrix index A}.  
\end{The}

\noindent {\bf Proof.} Let $U_i:=\Fbar_i(X_i)$, and $\Fbar_i(t)=t^{-\alpha_i}L_i(t)\in \mbox{RV}_{-\alpha_i}$, where $\lim_{t\to \infty}L_i(t)=l_i$, $1\le i\le d$. Since $\Fbar_i(t)/(t^{-\alpha_i}l_i) = L_i(t)/l_i\to 1$ as $t\to \infty$, it follows from Proposition 2.6 of \cite{Resnick07} that $\Fbar_i^{-1}(u)/(u^{-1/\alpha_i}l_i^{1/\alpha_i})\to 1$ as $u\to 0$. Note that Proposition 2.6 of \cite{Resnick07} is stated in terms of increasing regularly varying functions, but it can be easily verified that the result holds for decreasing regularly varying functions. That is, $\Fbar_i^{-1}(u)=u^{-1/\alpha_i}\ell_i(u)$ where $\ell_i(u)$ is slowly varying at $0$ satisfying that $\ell_i(u)\to l_i^{1/\alpha_i}$ as $u\to 0$.

For $u=t^{-1}$,  
consider, 
\by
\mathbb{P}(U_i\le u^{\alpha_i\gamma_i}x_i^{-\alpha_i}l_i, \ \exists i\in D) &=& \mathbb{P}(X_i> \Fbar_i^{-1}(u^{\alpha_i\gamma_i}x_i^{-\alpha_i}l_i), \ \exists i\in D) \\
&=& \mathbb{P}(X_i> (u^{\alpha_i\gamma_i}x_i^{-\alpha_i}l_i)^{-1/\alpha_i}\ell_i(u^{\alpha_i\gamma_i}x_i^{-\alpha_i}l_i), \ \exists i\in D) \\
&=& \mathbb{P}(X_i> t^{\gamma_i}x_il_i^{-1/\alpha_i}\ell_i(u^{\alpha_i\gamma_i}x_i^{-\alpha_i}l_i), \ \exists i\in D)
\ey
where $\lim_{u\to 0}l_i^{-1/\alpha_i}\ell_i(u^{\alpha_i\gamma_i}x_i^{-\alpha_i}l_i)= \lim_{u\to 0}l_i^{-1/\alpha_i}\ell_i(u^{\alpha_i\gamma_i})=1$, $1\le i\le d$. For any small $\epsilon>0$, we have that, as $u$ becomes sufficiently small, or equivalently, $t$ becomes sufficiently large, 
\[1-\epsilon\le l_i^{-1/\alpha_i}\ell_i(u^{\alpha_i\gamma_i}x_i^{-\alpha_i}l_i)\le 1+\epsilon,\ 1\le i\le d. 
\]
Hence, as $t$ is sufficiently large, or $u$ is sufficiently small,
\[\mathbb{P}(X_i> t^{\gamma_i}x_i(1+\epsilon), \ \exists i\in D)\le \mathbb{P}(X_i> t^{\gamma_i}x_il_i^{-1/\alpha_i}\ell_i(u^{\alpha_i\gamma_i}x_i^{-\alpha_i}l_i), \ \exists i\in D)
\]
\begin{equation}
\le \mathbb{P}(X_i> t^{\gamma_i}x_i(1-\epsilon), \ \exists i\in D).\label{TDF1}
\end{equation}
Dividing these two inequalities in \eqref{TDF1} by $R(t)=t^{-\beta}L(t)$ and noticing that $u=t^{-1}$, we have that for any small $\epsilon$, as $t$ is sufficiently large, or $u$ is sufficiently small,
\[\frac{\mathbb{P}(X_i> t^{\gamma_i}x_i(1+\epsilon), \ \exists i\in D)}{t^{-\beta}L(t)}\le \frac{\mathbb{P}(U_i\le u^{\alpha_i\gamma_i}x_i^{-\alpha_i}l_i, \ \exists i\in D)}{u^\beta L(u^{-1})}
\le \frac{\mathbb{P}(X_i> t^{\gamma_i}x_i(1-\epsilon), \ \exists i\in D)}{t^{-\beta}L(t)}.
\]
Taking limits as $t\to \infty$ and noticing that the intensity measure is continuous, we have
\[\lim_{u\to 0}\frac{\mathbb{P}(U_i\le u^{\alpha_i\gamma_i}x_i^{-\alpha_i}l_i, \ \exists i\in D)}{u^\beta L(u^{-1})}=\lim_{t\to \infty}\frac{\mathbb{P}(X_i> t^{\gamma_i}x_i, \ \exists i\in D)}{t^{-\beta}L(t)}= \mu\Big(\big(\prod_{i=1}^d[0,x_i]\big)^c\Big). 
\]
That is, $a_U(\ww; A, C)=a_L(\ww; A, \widehat{C})=\mu((\prod_{i=1}^d[0,w_i^{-1/\alpha_i}l_i^{1/\alpha_i}])^c)$ for any $\ww=(w_1, \dots, w_d)\in \mathbb{R}^d_+\backslash\{\boldsymbol{0}\}$. 
 \hfill
$\Box$

\medskip
If $L(t)\to r>0$, as $t\to \infty$ and if
\[\mu\Big([0, \infty]^{i-1}\times (1,\infty)\times [0,\infty]^{n-i}\Big)>0
\]  
then simple substitutions show that $\Fbar_i(t)\in \mbox{RV}_{-\alpha_i}$ where $\alpha_i=\beta/\gamma_i$, $1\le i\le d$. In this case, the matrix $A$ in \eqref{matrix index A} becomes the identity matrix $I$, and we have
\[a_U(\ww; I, C)=a_L(\ww; I, \widehat{C})=\mu\Big(\big(\prod_{i=1}^d[0,w_i^{-1/\alpha_i}l_i^{1/\alpha_i}]\big)^c\Big)
\]
for any $\ww=(w_1, \dots, w_d)\in \mathbb{R}^d_+\backslash\{\boldsymbol{0}\}$. That is, the copula $C$ has the tail dependence of order 1. 
Note that many multivariate distributions satisfy this property, as illustrated in the next example.

\begin{Exa}\rm 
	\label{4}
	Let $X_i=(T_i/Z)^{\gamma_i}$, $\gamma_i>0$, $i=1, 2$. where $(T_1, T_2)$ has a bivariate Marshall-Olkin distribution (see \eqref{MO survival})
	\[\mathbb{P}(T_1>t_1, T_2>t_2)= \exp\{-t_1-t_2-\lambda\max\{t_1,t_2\}\},\ \ \lambda>0,
	\]
	and $Z$, independent of $(T_1, T_2)$, has a gamma distribution with scale parameter 1 and shape parameter $\beta>0$. The distribution of $(X_1, X_2)$ is known as a bivariate Pareto of the fourth kind (see page 586, \cite{KBJ00}) and its survival function is given by
	\begin{equation}
	\label{Pareto}
	\overline{F}(x_1, x_2)=\left[1+x_1^{1/\gamma_1}+x_2^{1/\gamma_2}+\lambda\max\{x_1^{1/\gamma_1}, x_2^{1/\gamma_2}\}\right]^{-\beta},\ x_1\ge 0, x_2\ge 0. 
	\end{equation}
	The marginal survival functions are given by
	\[\overline{F}_i(t)=\left[1+(1+\lambda)t^{1/\gamma_i}\right]^{-\beta}, \ t\ge 0, \ i=1, 2. 
	\]
	Plug the inverse $\overline{F}_i^{-1}(u)$, $i=1,2$, into \eqref{Pareto}, and we obtain the survival copula
	\begin{equation}
	\label{Pareto copula}
	\widehat{C}(u_1, u_2)=
	\left[1+\Big(\frac{u_1^{-1/\beta}-1}{1+\lambda}\Big)+\Big(\frac{u_2^{-1/\beta}-1}{1+\lambda}\Big)+\lambda\max\Big\{\frac{u_1^{-1/\beta}-1}{1+\lambda}, \frac{u_2^{-1/\beta}-1}{1+\lambda}\Big\}\right]^{-\beta}.
	\end{equation}
	To derive the intensity measure for $(X_1, X_2)$, consider
	\[
	\begin{aligned}
	&\mathbb{P}(X_1>t^{\gamma_1}x_1\ \mbox{or}\  X_2>t^{\gamma_2}x_2)=\sum_{i=1}^2\mathbb{P}(X_i>t^{\gamma_i}x_i)-\overline{F}(t^{\gamma_1}x_1,t^{\gamma_2}x_2)\\
	&= \sum_{i=1}^2\left[1+t(1+\lambda)x_i^{1/\gamma_i}\right]^{-\beta}
	-\left[1+tx_1^{1/\gamma_1}+tx_2^{1/\gamma_2}+t\lambda\max\{x_1^{1/\gamma_1}, x_2^{1/\gamma_2}\}\right]^{-\beta}.
	\end{aligned}
	\]
	The intensity measure is given by
	\[
		\begin{aligned}
		& \mu\big(([0,x_1]\times [0,x_2])^c\big) = 
	\lim_{t\to \infty}\frac{\mathbb{P}(X_1>t^{\gamma_1}x_1\ \mbox{or}\  X_2>t^{\gamma_2}x_2)}{t^{-\beta}}\\
	&=\sum_{i=1}^2\left[(1+\lambda)x_i^{1/\gamma_i}\right]^{-\beta}
	-\left[x_1^{1/\gamma_1}+x_2^{1/\gamma_2}+\lambda\max\{x_1^{1/\gamma_1}, x_2^{1/\gamma_2}\}\right]^{-\beta}
		\end{aligned}
	\]
	and the scaling function $R(t)=t^{-\beta}$ and the marginal scaling indexes are $\gamma_1$ and $\gamma_2$. Note that 
	\[\overline{F}_i(t)=t^{-\beta/\gamma_i}L_i(t), \ \mbox{where}\ L_i(t)=\frac{\left[1+(1+\lambda)t^{1/\gamma_i}\right]^{-\beta}}{t^{-\beta/\gamma_i}}\to (1+\lambda)^{-\beta}, \ \mbox{as}\ t\to \infty. 
	\]
	On the other hand, it follows from \eqref{Pareto copula} that
	\[
	\begin{aligned}
		& a_L(\ww; I, \widehat{C})=w_1+w_2-\lim_{u\to 0}\frac{\widehat{C}(uw_1, uw_2)}{u}\\
		& =  w_1+w_2-\left[\frac{1}{1+\lambda}(w_1^{-1/\beta}+w_2^{-1/\beta})+\frac{\lambda}{1+\lambda}\max\{w_1^{-1/\beta}, w_2^{-1/\beta}\}\right]^{-\beta}.
	\end{aligned}
	\]
	It is easy to verify that 
	\[a_L(\ww; I, \widehat{C})= \mu\Big(\big([0,w_1^{-\gamma_1/\beta}(1+\lambda)^{-\gamma_1}]\times [0,w_1^{-\gamma_2/\beta}(1+\lambda)^{-\gamma_2}]\big)^c\Big), 
	\]
	for any $w_1>0, w_2>0$. 
	\hfill
	$\Box$
\end{Exa} 

It is worth mentioning that if for some $i$ in Theorem \ref{MRV-TDF}
\[\mu\Big([0, \infty]^{i-1}\times (1,\infty)\times [0,\infty]^{n-i}\Big)=0
\]  
then the $i$-th marginal tail is lighter than other marginal tails and the exponent function $a_U(\cdot; A, C)$ describes the dependence hidden within the subcone $\mathbb{E}^{(2)}$. 

\section{Concluding Remarks}

We introduce in this paper the operator tail dependence for copulas, and show that non-standard multivariate regular variation can be characterized by copulas with operator tail dependence.  
The advantage of the copula method is due to separation of univariate margins from dependence analysis, and this is illustrated by Theorems \ref{TDF-MRV} and \ref{MRV-TDF} in the context of multivariate extreme value analysis. For non-standard multivariate regular variation, we can use copulas with operator tail dependence, incorporated with regularly varying univariate margins with possibly different tail indexes. 

Operator regular variation has been motivated by the observations that multivariate data with heavy tails may have different tail indexes along different directions (see \cite{MS01, BE07}), where the directions
are not necessarily along the original coordinate axes. A statistical test was developed in \cite{MS2013} to determine whether the tail index varies with direction in any given data set, so that operator norming can be used. As illustrated in \cite{NPM01, MS99}, pooling tail estimates from each original marginal variable can be misleading, since it fails to detect the lighter tails along certain direction. Operator regular variation, such as non-standard multivariate regular variation in possibly rotated coordinate systems, can determine completely tail behaviors. 
It must be mentioned that copulas may not be invariant under orthogonal transforms. For any positive-definite matrix $A$ of tail indexes (see \eqref{dec}), we need to take the orthogonal transform (such as a rotation) on the original vector data and then use the copula method on the transformed vector.

\bigskip
\noindent
{\bf Acknowledgments:} The author would like to sincerely thank referees and an associate editor for their insightful comments, which led to an improvement of the presentation and motivation of this paper. The author would also like to thank Mark Meerschaert for a useful discussion.


\end{document}